\documentclass[12pt]{amsart}
\usepackage{lipsum}
\usepackage{amssymb}
\usepackage{enumerate}
\usepackage{amsthm}
\usepackage{tikz-cd}
\usepackage{amsmath}
\usepackage[mathscr]{euscript}
\usepackage[utf8]{inputenc}
\usepackage[english]{babel}
\usepackage{hyperref}
\hypersetup{
    colorlinks=true,
    linkcolor=blue,
    citecolor=red,
    filecolor=red,      
    urlcolor=violet,
}

\makeatletter
\@namedef{subjclassname@2020}{%
  \textup{2020} Mathematics Subject Classification}
\makeatother
\newtheorem{thm}{Theorem}[section]
\newtheorem{lem}[thm]{Lemma}
\newtheorem{prop}[thm]{Proposition}

\newtheorem{corl}[thm]{Corollary}
\newtheorem{xrem}{Remark}
\newtheorem{exm}[thm]{Example}
\newtheorem{conj}[thm]{Conjecture}
\newtheorem{prob}[thm]{Problem}

\DeclareMathOperator{\Pic}{{Pic}}

\DeclareMathOperator{\Div}{{Div}}

\DeclareMathOperator{\Bl}{{Bl}}
\DeclareMathOperator{\Eff}{{Eff}}

\DeclareMathOperator{\Sym}{{Sym}}

\begin{document}
\baselineskip=14pt
\author{Snehajit Misra}
\author{Nabanita Ray}
\address{Indraprastha Institute of Information Technology, Delhi (IIITD),
Okhla Phase-III, New Delhi-110020, India.}
\email[Snehajit Misra]{misra08@gmail.com}
\address{Indraprastha Institute of Information Technology, Delhi (IIITD),
Okhla Phase-III, New Delhi-110020, India.}
\email[Nabanita Ray]{nabanitaray2910@gmail.com}

\begin{abstract}
In the first part of this article, we give bounds on self-intersections $C^2$ of integral curves $C$ on  blow-ups $\Bl_nX$ of surfaces $X$ with the anti-cannonical divisor $-K_X$ effective. In the last part, we prove the weak bounded negativity for self-intersections  $C^2$ of integral curves $C$ in a family of surfaces $f:Y\longrightarrow B$ where $B$ is a smooth curve.
\end{abstract}

\title{On Weak bounded negativity conjecture}
\maketitle

\section{Introduction}
The so-called weak bounded negativity conjecture (see \cite[Theorem 2.1]{H19}) is derived from the study of open problem named ``bounded negativity conjecture (BNC)" which is stated as follows :
\begin{conj}
 $($Bounded Negativity Conjecture$)$ For every smooth projective surface $X$ over the complex numbers, there exists a non-negative integer $b(X)\in \mathbb{Z}$ such that $$C^2\geq - b(X)$$
 for all integral curves $C\subset X$.
\end{conj}

Moreover, we have a stronger bounded negativity conjecture for the blow-up surface $\Bl_{n}\mathbb{P}^2$  of $\mathbb{P}^2$ at $n$ very general points.
\begin{conj}\label{conj1.2}
 If $C$ is an integral curve in $\Bl_{n}\mathbb{P}^2$, then
 \begin{center}
 $C^2\geq -1 $ and  $C^2 = -1$ implies $C$ is a (-1)-curve.
 \end{center}
 \end{conj}
 Another motivation for studying this stronger version is the following Nagata conjecture. The stronger form of bounded negativity conjecture for $\Bl_{n}\mathbb{P}^2$ implies Nagata conjecture:
 \begin{conj}
  $($Nagata conjecture$)$ Let $p_1,p_2,\cdots,p_n\in \mathbb{P}^2$ be $n$ very general points with $n\geq 9$. If $C\subset \mathbb{P}^2$ is an integral curve with multiplicity $m_i$ at $p_i$, then the degree of $C$ satisfies the following :
  $$\deg(C)\geq \frac{1}{\sqrt{n}}\sum\limits_{i=1}^nm_i.$$
 \end{conj}
 Also SHGH Conjecture implies the Conjecture \ref{conj1.2} (see \cite{C13}). Hence it is worth studying the bounded negativity conjecture for rational
surfaces.

Another interesting problem posed by Demailly asks whether the global Seshadri constant is positive for any smooth surface $X$ (see \cite[Question 6.9]{D92}). We denote the Seshadri constant of a line bundle $L$ on a smooth surface $X$ at a point $x\in X$ by $\varepsilon(X,L;x)$. It is still unknown if for every fixed $x\in X$, the value
\begin{center}
$\inf\bigl\{\varepsilon(X,L;x) \mid L\in \Pic(X)$ is ample $\bigr\}$
\end{center} is always positive. Infact, BNC implies an affirmative answer to this question (see \cite[Section 3.6]{BBCRDHJKKRR12}).

 The following weaker version of Bounded Negativity Conjecture was proposed in \cite[Conjecture 3.3.4]{BBCRDHJKKRR12}.

 \begin{thm} \cite[Theorem 2.1]{H19}
 $($Weak Bounded Negativity conjecture$)$ For any smooth complex projective surface $X$ and any integer $g$, there is a constant $b(X,g)$, only depending on $X$ and the integer $g$, such that
 $$C^2\geq b(X,g)$$ for any reduced curve $C=\sum\limits_i C_i$ in $X$ with the geometric genus $p_g(C_i)\leq g$ for all $i$.
\end{thm}

The authors in \cite{BBCRDHJKKRR12} used the logarithmic
Miyaoka-Yau inequality to prove the weak bounded negativity conjecture for smooth surfaces with non-negative Kodaira dimension and gave the following bound.
\begin{thm}\cite[Theorem 2.6]{BBCRDHJKKRR12}
 Let $X$ be a smooth projective surface with Kodaira dimension $\kappa(X)\geq 0$. Then for any integral curve $C\subset X$, we have
 $$C^2\geq K_X^2-3c_2(X)+2-2p_g(C)$$
 where $c_2$ is the second Chern number of the surface $X$ and $p_g(C)$ is the geometric genus of $C$.
\end{thm}

Subsequently, the weak bounded negativity conjecture is proved  for all smooth projective surfaces (see \cite[Corollary 1.12, Theorem 2.1]{H19}). More specifically, the author in \cite{H19} proved the following Theorem \ref{thm4} by making use of logarithmic Miyaoka-Yau inequality and Zariski decomposition of pseudoeffective divisors on smooth projective surfaces.

\begin{thm}\cite[Corollary 1.12]{H19}\label{thm4}
  Let $X$ be a smooth projective surface  and let $C$ be an integral curve of geometric genus $p_g(C)$ on $X$. Consider the following cases :
  \begin{itemize}
   \item \bf Case 1 : \it Suppose $h^0(-K_X) \neq 0$. Then, in this case, the self-intersection $C^2\geq -2$ for all integral curves $C$ except finitely many components of the linear system $\lvert -K_X\rvert.$
   \vspace{2mm}

   \item \bf Case 2 : \it Suppose $h^0(-K_X) = 0$. Then we consider the following two subcases.
   \vspace{1mm}

   \begin{enumerate}
    \item \bf Subcase 2.1 : \it Suppose $C$ is an integral curve such that $h^0\bigl(2(K_X+C)\bigr) = 0$. Then $$C^2\geq K_X^2+\chi(\mathcal{O}_X)-3.$$
    \item \bf Subcase 2.2 : \it Suppose $C$ is an integral curve such that $h^0\bigl(2(K_X+C)\bigr) \neq  0$. Then
    $$C^2\geq K_X^2-3c_2(X)+2-2p_g(C).$$
   \end{enumerate}
\end{itemize}
Combining Case 1 and Case 2, the weak bounded negativity conjecture is proved for all surfaces.
 \end{thm}

In the recent past, there are works done for the following weighted version of the BNC, which is usually named the  weighted BNC, which we will state now.

\begin{conj}\label{conj1.7}
 $($Weighted Bounded Negativity Conjecture$)$ For every smooth projective surface $X$ over the complex numbers, there exists a non-negative integer $b_{\omega}(X) \in \mathbb{Z}$ such that $$C^2\geq -b_{\omega}(X)(C\cdot H)^2$$ for all integral curves $C \subset X$ and all big and nef line bundles $H$ for which $C\cdot H >0$.
\end{conj}
The previous Conjecture \ref{conj1.7} asserts that the weighted self-intersection $C^2/(C\cdot H)^2$ of $C$ to be bounded from below. Putting in a different perspective, the above conjecture asking for a bound on the self-intersection  $C^2$ of all integral curves $C$ on $X$ that depends on both $X$ and the degree of the curve $C$ with respect to every big and nef line bundle over which the integral curve is positive. The authors in \cite{LP20} has obtained the following result regarding the weighted bounded negativity conjecture for surfaces with non-negative Kodaira dimension.
\begin{thm}\cite[Theorem A]{LP20}
 Assume that $X$ is a surface of non-negative Kodaira dimension and let $f:Y\longrightarrow X$ be the blowing up of $X$ along $n$ mutually distinct points. Then there exists a big and nef line bundle $\Gamma$ such that
 $$C^2\geq -\frac{1}{2}\bigl(\delta(X)+C\cdot \Gamma)-n,$$
 for every integral curve $C\subset Y$, where $\delta(X) = 3e(X)-K_X^2.$
 \end{thm}
The authors in \cite{LP20} have also obtained certain bounds for the self-intersections of integral curves in blow-up rational surface $\Bl_n\mathbb{P}^2$ and $\Bl_n\mathbb{F}_e$, where $\mathbb{F}_e$ is the Hirzebruch surfaces with invariants $e$ (see \cite[Theorem B]{LP20}).

Subsequently, the bounds for self-intersection of integral curves in $\Bl_n\mathbb{P}^2$ are improved in \cite{CF24}.  The authors in \cite{CF24} proved the following.
\begin{thm}\cite[Theorem 2]{CF24}
 Let $Y$ be the blow-up of $\mathbb{P}^2$ at $n$-arbitrary (proper or infinitely near) points and let $L$ be the pullback of a (general) line in $\mathbb{P}^2$. Then all but finitely many integral curves $C$ on $Y$ satisfy the inequality $$C^2\geq \min\Bigl\{-2,\hspace{1mm} -\frac{1}{12}n(C\cdot L+27)\Bigr\}.$$
\end{thm}

In the present paper, we focus on integral curves on blow-ups of certain surfaces and bounding their negativity. Our main results are as follows:
\begin{thm}\label{thm1.11}
 Let $X$ be a smooth projective irreducible surface with $\chi(\mathcal{O}_X) \geq 1$ and  $-K_X$ is effective. Furthermore, let $L$ be a very ample divisor on $X$ with the intersection number $L\cdot K_X = - a_0 < 0$. Let $\pi_n : X_n = \Bl_n X \longrightarrow X$ be the blow-up map at $n$ points, and $H=\pi_n^*L$. Then for any integral curve $C$ on $X_n$, we have the following two cases :
 \begin{enumerate}
  \item \bf Case 1 : \it If $K_X^2\leq n$, then $$C^2 \geq \min\bigl\{\mathcal{L},\mathcal{M}\bigr\},$$
  \vspace{2mm}

 \item \bf Case 2 : \it If $K_X^2> n$, then
$$C^2 \geq \min\bigl\{\mathcal{M},\mathcal{N}\bigr\},$$
where
$$\mathcal{L}= \frac{(C\cdot H+a_0)}{2a_0}\Bigl(K_{X}^2-n\Bigr)-3, $$
$$\mathcal{N}= \frac{(C\cdot H+1)}{2a_0}\Bigl(K_{X}^2-n\Bigr)-3,$$
  and $$\mathcal{M}= \frac{1}{2}(H^2+1)\bigl(K_{X}^2 -n\bigr) - a_0^2 -3 + \frac{a_0C\cdot H}{H^2}. $$
\end{enumerate}
\end{thm}

\begin{thm}\label{thm1.12}
 Let $X$ be a smooth projective irreducible surface with $\chi(\mathcal{O}_X) < 1$  and  $-K_X$ is effective. Furthermore, let $L$ be a very ample divisor on $X$ with the intersection number $L\cdot K_X = - a_0 < 0$. Let $\pi_n : X_n = \Bl_n X \longrightarrow X$ be the blow-up map at $n$ points, and $H=\pi_n^*L$. Then for any integral curve $C$ on $X_n$, we have the following two cases :
 \begin{enumerate}
  \item \bf Case 1 : \it If $K_X^2\leq n$, then $$C^2 \geq \min\bigl\{\mathcal{L},\mathcal{M}\bigr\},$$
  \vspace{2mm}

  \item \bf Case 2 : \it If $K_X^2> n$, then
$$C^2 \geq \min\bigl\{\mathcal{M},\mathcal{N}\bigr\},$$
where
$$\mathcal{L}= \chi(\mathcal{O}_X)+\frac{(C\cdot H +a_0)}{2a_0}\Bigl(K_X^2-n\Bigr)-4, $$
$$\mathcal{N}=  \chi(\mathcal{O}_X)  + \frac{(C\cdot H+1)}{2a_0}\Bigl(K_X^2-n\Bigr) -4,$$
  and $$\mathcal{M}= \frac{1}{2}(H^2+1)\bigl(K_{X}^2 -n\bigr) - a_0^2 -3 + \frac{a_0C\cdot H+\chi(\mathcal{O}_X)-1}{H^2}.$$
  \vspace{1mm}

 \end{enumerate}
 \end{thm}
As  applications of our results, we obtain better bounds in the case of $\Bl_n\mathbb{P}^2$, $\Bl_n\mathbb{F}_e$ and blow-up of certain ruled surfaces (see Example \ref{exm3.6} - Example \ref{exm3.8}). We observe that the results in these Theorem \ref{thm1.11} and Theorem \ref{thm1.12} can be proved in a more general set up whenever some integral multiple $-mK_X$ of the anti-cannonical divisor $-K_X$ is effective. In other words, the ideas in Theorem \ref{thm1.11} and Theorem \ref{thm1.12} holds true for the base surfaces $X$ having Kodaira dimension $-\infty$, which complements the results in \cite{LP20}.

Next, in the last section, we consider the following question which is  asked in \cite[Problem 3.3.6]{BBCRDHJKKRR12}.
\begin{prob}\label{prob1}
 Let $f : Y \longrightarrow B$ be a morphism from a smooth projective threefold $Y$ to a smooth curve $B$ such that the general fibre is a smooth surface. Is there a
constant $b(Y,g)$ such that
$$C^2\geq-b(Y,g)$$
for all vertical curves $C \subset Y$ (i.e., $f(C)$ = a point) of geometric genus $g$? (Here the
self-intersection is computed in the fibre of $f$ containing $C$.)
\end{prob}
This Problem \ref{prob1} motivates us to prove the following.
 \begin{thm}
Let $f:Y\rightarrow B$ be a morphism from a smooth projective threefold $Y$ to a smooth curve $B$ such that generic  fibers are smooth surfaces. Let $C=\sum\limits_{i=1}^m C_i\subset Y$ be any reduced vertical curve such that $C$ is contained in a smooth fiber $Y_x$. Let $g=\max\bigl\{p_g(C_i) \mid 1\leq i \leq m\bigr\}$. Then there exists a constant $B(Y,g)$ only depending on $Y$ and
$g$, such that
$$C^2 \geq  B(X,g)$$.
\end{thm}

\section{Notations and Preliminaries}
Throughout this article, all the algebraic varieties are assumed to be irreducible and defined over the field of complex numbers $\mathbb{C}$. Given a divisor $D$ on a smooth projective variety $X$, we write $h^i(X,D)$ to denote the dimension of $i$-th cohomology group $H^i\bigl(X,\mathcal{O}_X(D)\bigr)$ as a complex vector space. We will simply write $h^i(D)$ whenever the base variety $X$ is clear from the context. The sheaf $K_X$ will denote the canonical sheaf on $X$. The arithmetic and geometric genus of an integral curve $C$ will be denoted by $p_a(C)$ and $p_g(C)$ respectively. The dual of a vector bundle $V$ will be denoted by $V^{\vee}$.

Let $X$ be a smooth projective variety over $\mathbb{C}$ and the set $\Div(X)$ denotes the set of all divisors on $X$. Let
  \begin{center}
  $\Div^0(X) := \bigl\{ D \in \Div(X) \mid D\cdot C = 0 $ for all curves $C$ in $X \bigr\} \subseteq \Div(X)$.
  \end{center}
 be the subgroup of $\Div(X)$ consisting of numerically trivial divisors. The quotient $\Div(X)/\Div^0(X)$ is called the N\'{e}ron Severi group of $X$, and is denoted by $N^1(X)_{\mathbb{Z}}$.
   The N\'{e}ron Severi group  $N^1(X)_{\mathbb{Z}}$ is a free abelian group of finite rank.
 Its rank, denoted by $\rho(X)$ is called the Picard number of $X$. In particular, $N^1(X)_{\mathbb{R}}$ is called the real N\'{e}ron
 Severi group and $N^1(X)_{\mathbb{R}}  := N^1(X)_{\mathbb{Z}} \otimes \mathbb{R} := \bigl(\Div(X)/\Div^0(X)\bigr) \otimes \mathbb{R}$.

 The convex cone generated by the set of all effective divisors in $N^1(X)_\mathbb{R}$ is denoted by $\Eff^1(X)$ and its closure $\overline{\Eff}^1(X)$ is called the \it pseudo-effective cone \rm  of divisors in $X$. The elements of $\overline{\Eff}^1(X)$ are called pseudo-effective divisors on $X$.
\subsection{Zariski-decomposition of pseudo-effective divisors} Let $X$ be a smooth projective surface, and let $D$ be a pseudo-eﬀective integral divisor on $X$. Then $D$ can be written
uniquely as a sum $D = P + N$ of $\mathbb{Q}$-divisors with the following properties:
\begin{enumerate}
 \item $P$ is nef;
 \item $N = \sum\limits_{i=1}^m a_iE_i$ is effective, and if $N\neq 0$, then the intersection matrix $(E_i\cdot E_j)_{1\leq i,j\leq m}$ determined by the components of $N$ is negative deﬁnite;
 \item $P$ is orthogonal to each of the components of $N$ , i.e., $P \cdot E_i = 0$ for each $i$.
\end{enumerate}
In the following statement, $\left\lceil mN\right\rceil$ denotes the round-up of $mN$ , i.e. the divisor
obtained from $mN$ by rounding up the coefficients of each of its components.
In the situation of a Zariski decomposition as above, the natural map $$H^0\bigl(X,
\mathcal{O}_X(mD-\left\lceil mN\right\rceil)\bigr)\longrightarrow H^0(X,\mathcal{O}_X(mD))$$
is bijective for every integer $m\geq 1$.
\section{Bounding negativity of integral curves on certain blow-ups}
We begin this section by the following easy lemma which is proved in \cite[Lemma 6]{CF24}. For the sake of completeness we give the proof of the following lemma.
\begin{lem}\label{lem2.1}
 Let $X$ be a smooth projective surface and $C$ be an integral curve on $X$. Then for any integer $m\neq 1$, we have
\begin{align*}
   C^2
  &= \frac{\chi(\mathcal{O}_X)}{m-1}+\frac{1}{2}mK_X^2+2p_a(C)+
 \frac{p_a(C)}{m-1}-2-
 \frac{1}{m-1}-\frac{1}{m-1}h^0(mK_X+C) \\
&\hspace{2cm} + \frac{1}{m-1}h^1(mK_X+C)-\frac{1}{m-1}h^0\bigl(-(m-1)K_X-C\bigr).
\end{align*}
\begin{proof}
By Riemann-Roch Theorem for the divisor  $mK_X+C$ on smooth surface $X$ we have
\begin{align*}
& h^0(mK_X+C) - h^1(mK_X+C) + h^2(mK_X+C) \\
& h^0(mK_X+C) - h^1(mK_X+C) + h^0(-(m-1)K_X-C)\\
&= \frac{1}{2}\Bigl\{(mK_X+C)\bigl((m-1)K_X+C\bigr)\Bigr\} + \chi(\mathcal{O}_X).
\end{align*}
By adjunction formula we have $$C^2=2p_a(C)-2-K_X\cdot C.$$
Thus
\begin{align*}
 & \frac{1}{2}\Bigl\{(mK_X+C)\bigl((m-1)K_X+C\bigr)\Bigr\}\\
 & = \frac{1}{2}\Bigl[m(m-1)K_X^2+(2m-1)K_X\cdot C+C^2\Bigr]\\
 & = \frac{1}{2}\Bigl[m(m-1)K_X^2+(2m-1)(2p_a(C)-2-C^2)+C^2\Bigr]\\
 &= \frac{1}{2}\Bigl[m(m-1)K_X^2+(2m-1)(2p_a(C)-2)-C^2(2m-2) \Bigr]\\
 &= (m-1)\Bigl[\frac{1}{2}mK_X^2+\frac{(2m-1)(p_a(C)-1)}{(m-1)} - C^2\Bigr]\\
 & = (m-1)\Bigl[\frac{1}{2}mK_X^2+\frac{\bigl\{2(m-1)+1\bigr\}(p_a(C)-1)}{(m-1)} - C^2\Bigr]\\
 &=(m-1)\Bigl[\frac{1}{2}mK_X^2+\bigl\{2+\frac{1}{m-1}\bigr\}(p_a(C)-1) - C^2\Bigr]\\
 &=(m-1)\Bigl[\frac{1}{2}mK_X^2+ 2p_a(C)-2+\frac{1}{m-1}p_a(C)-\frac{1}{m-1} - C^2 \Bigr]
\end{align*}
Therefore we conclude that
\begin{align*}
 C^2
  &= \frac{\chi(\mathcal{O}_X)}{m-1}+\frac{1}{2}mK_X^2+2p_a(C)+
 \frac{p_a(C)}{m-1}-2-
 \frac{1}{m-1}-\frac{1}{m-1}h^0(mK_X+C) \\
& \hspace{2cm} + \frac{1}{m-1}h^1(mK_X+C)-\frac{1}{m-1}h^0\bigl(-(m-1)K_X-C\bigr).
\end{align*}
\end{proof}
\end{lem}
We use the following  notations throughout the rest of this section. Let $X$ be a smooth projective surface, $\pi_n:X_n\longrightarrow X$ be the blow-up of $X$ at $n$ distinct points, and $L$ be a very ample divisor on $X$. Note that $\pi_n^*L$ is a big and nef divisor on $X_n$.  We denote $\pi_n^*L$ as $H$ on $X_n$, and the exceptional divisors of the map $\pi_n$ as $E_1,E_2,\cdots, E_n$.

\begin{xrem}\label{xrem3.2}
    Let $E$ and $E'$ be two  exceptional divisors on $X_n$ which are not necessarily distinct. Then observe that $\dim_{\mathbb{C}} H^0(X_n,\mathcal{O}_{X_n}(E\otimes E'))=1$. This follows from the exact sequence
    $$0\longrightarrow\mathcal{O}_{X_n}(E')\longrightarrow\mathcal{O}_{X_n}(E\otimes E')\longrightarrow\mathcal{O}_E(-1)=\mathcal{O}_{\mathbb{P}^1}(-1)\longrightarrow 0$$
\end{xrem}
\begin{lem}\label{lem3.3}
Let $D$ be any effective divisor on $X_n$ such that $D\cdot H=0$. Then $h^0(X_n,D)\leq 1$.
\begin{proof}
    Any effective divisor $D$ on $X_n$ is linearly equivalent to $F+\sum a_iE_i$, where $F$ is the pull back of some effective divisor from $X$ and $a_i\in \mathbb{Z}$. As $H\cdot D=0$, $F=0$. Hence $D$ is linearly equivalent to linear combination of some exceptional divisors of the map $\pi_n$. Therefore, from Remark \ref{xrem3.2}, $h^0(X_n,D)\leq 1$.
\end{proof}
\end{lem}

\begin{lem}\label{lem2.2}
Let $X$ be a smooth projective irreducible surface such that $-K_X$ is effective, and $L$ be a very ample divisor on $X$ with the intersection number $L\cdot K_X = - a_0 < 0$. We denote the blow-up of the surface $X$ at $n$ points by $X_n$, and let $C$ be any curve  (not necessarily integral) on $X_n$. Then there exists a positive integer $m_C$  which satisfies the following inequality
 $$\frac{C\cdot H+1}{a_0}\leq m_C\leq \frac{C\cdot H + a_0}{a_0}.$$
 Furthermore $$h^0(m_CK_{X_n}+C) = 0 \hspace{2mm} \text{and} \hspace{2mm} 0\leq  h^2 (m_CK_{X_n}+C) = h^0\bigl(-(m_C-1)K_{X_n}-C\bigr)  \leq 1.$$

 \begin{proof}
Let $m_C$ be the smallest positive integer such that
\begin{align}\label{seq 1}
(m_CK_{X_n}+C)\cdot H \leq -1.
\end{align}
i.e.
$$
C\cdot H -a_0m_C \leq -1.
$$

Equivalently,
\begin{align}\label{seq 2}
 \frac{C\cdot H +1}{a_0} \leq m_C.
\end{align}
Also by minimality of $m_C$, we have
$$\bigl((m_C-1)K_{X_n}+C\bigr)\cdot H \geq 0.$$
Therefore
 $$-a_0\leq C\cdot H -a_0m_C.$$
 Thus
\begin{align}\label{seq 3}
m_C\leq \frac{C\cdot H+a_0}{a_0}.
\end{align}
Combining ($\ref{seq 2}$) and ($\ref{seq 3}$) we get the following inequalities:
$$\frac{C\cdot H+1}{a_0}\leq m_C\leq \frac{C\cdot H + a_0}{a_0}.$$

From ($\ref{seq 1}$) we also conclude that $m_CK_{X_n}+C$ is not linearly equivalent to an effective divisor.

Hence $$h^0\bigl(m_CK_{X_n}+C\bigr) = 0.$$
Now by the minimality of $m_C$, we deduce from ($\ref{seq 1}$) that $$\bigl((m_C-1)K_{X_n}+C\bigr)\cdot H \geq 0.$$
Therefore using Lemma \ref{lem3.3} we have
 $$0\leq h^0\bigl(-(m_C-1)K_{X_n}-C\bigr) = h^2\bigl(m_CK_{X_n}+C\bigr) \leq 1.$$
This completes the proof.
\end{proof}
\end{lem}

\begin{thm}\label{thm3.2}
 Let $X$ be a smooth projective irreducible surface with $\chi(\mathcal{O}_X) \geq 1$ and  $-K_X$ is effective. Furthermore, let $L$ be a very ample divisor on $X$ with the intersection number $L\cdot K_X = - a_0 < 0$. Let $\pi_n : X_n = \Bl_n X \longrightarrow X$ be the blow-up map at $n$ points, and $H=\pi_n^*L$. Then for any integral curve $C$ on $X_n$, we have the following two cases :
 \begin{enumerate}
  \item \bf Case 1 : \it If $K_X^2\leq n$, then $$C^2 \geq \min\bigl\{\mathcal{L},\mathcal{M}\bigr\},$$
  \vspace{2mm}

 \item \bf Case 2 : \it If $K_X^2> n$, then
$$C^2 \geq \min\bigl\{\mathcal{M},\mathcal{N}\bigr\},$$
where
$$\mathcal{L}= \frac{(C\cdot H+a_0)}{2a_0}\Bigl(K_{X}^2-n\Bigr)-3, $$
$$\mathcal{N}= \frac{(C\cdot H+1)}{2a_0}\Bigl(K_{X}^2-n\Bigr)-3,$$
  and $$\mathcal{M}= \frac{1}{2}(H^2+1)\bigl(K_{X}^2 -n\bigr) - a_0^2 -3 + \frac{a_0C\cdot H}{H^2}. $$
\end{enumerate}
 \vspace{2mm}

 \begin{proof}
   We consider the following two cases depending on whether an integral curve $C$ has $m_C =1$  or not (see Lemma \ref{lem2.2} for the notation $m_C$) :
  \begin{enumerate}
   \item \bf Case 1 : \rm Let $C$ be an integral curve on $X$ such that $m_C \neq 1$. Then using Lemma \ref{lem2.1} we get
  \begin{align*}
   & C^2 = \frac{\chi(\mathcal{O}_X)}{m_C-1}+\frac{1}{2}m_CK_X^2+2p_a(C)+
 \frac{p_a(C)}{m_C-1}-2-
 \frac{1}{m_C-1} \\
  & \hspace{1cm} -\frac{1}{m_C-1}h^0(m_CK_X+C) +
 \frac{1}{m_C-1}h^1(m_CK_X+C)-\frac{1}{m_C-1}h^0\bigl(-(m_C-1)K_X-C\bigr).
\end{align*}
Note that by Lemma $\ref{lem2.2}$ we have $$h^0(m_CK_X+C)=0 \hspace{2mm} \text{and} \hspace{2mm}
h^0\bigl(-(m_C-1)K_X-C)\leq 1.$$
This implies
\begin{align}\label{seq1}
C^2 \geq \frac{\chi(\mathcal{O}_X)-1}{m_C-1}+\frac{1}{2}m_CK_{X_n}^2-2+2p_a(C)-1.
\end{align}
By the hypothesis we have $\chi(\mathcal{O}_X) \geq 1$. Thus
\begin{align}\label{s1}
C^2 \geq \frac{1}{2}m_CK_{X_n}^2-3 = \frac{1}{2}m_C(K_{X}^2-n)-3.
\end{align}
\vspace{1mm}

We now consider the following two subcases :
\vspace{2mm}

\begin{enumerate}
 \item \bf Subcase 1 : \rm Suppose that $K_X^2\leq n$. Recall that from Lemma \ref{lem2.2} $$m_C\leq \frac{C\cdot H+a_0}{a_0}.$$ Thus, in this case, using $(\ref{s1})$ we get
 $$C^2\geq \frac{(C\cdot H+a_0)}{2a_0}\Bigl(K_{X}^2-n\Bigr)-3.$$
 \vspace{1mm}

 \item \bf Subcase 2 : \rm Suppose that $K_X^2 > n$. Then again using Lemma \ref{lem2.2} we have $$m_C\geq \frac{C\cdot H+1}{a_0}.$$
 Therefore in this case we get using $(\ref{s1})$ $$C^2\geq \frac{(C\cdot H+1)}{2a_0}\Bigl(K_{X}^2-n\Bigr)-3.$$
\end{enumerate}
\vspace{3mm}

\item \bf Case 2 : \rm Let $C$ be an integral curve on $X_n$  such that $m_C=1$. Let us fix the reduced curve $C'= C+ a_0H$. Again using Lemma \ref{lem2.2} we have the following inequality:
$$\frac{C'\cdot H+1}{a_0} \leq m_{C'}\leq \frac{C'\cdot H+a_0}{a_0}.$$
This implies
$$\frac{\Bigl\{(C+a_0H)\cdot H\Bigr\}+1}{a_0} \leq m_{C'}\leq \frac{\Bigl\{(C+a_0H)\cdot H\Bigr\}+a_0}{a_0}.$$
i.e.
\begin{align}\label{seq 4}
H^2 + \frac{C\cdot H+1}{a_0}\leq m_{C'} \leq  H^2 + \frac{C\cdot H+a_0}{a_0}
\end{align}
We also have $$\frac{C\cdot H+1}{a_0}\leq m_C = 1 \leq \frac{C\cdot H + a_0}{a_0}.$$
Therefore
\begin{align}\label{seq 5}
-\frac{C\cdot H + a_0}{a_0}\leq - 1 \leq -\frac{C\cdot H+1}{a_0}
\end{align}
Adding $(\ref{seq 4})$ and $(\ref{seq 5})$ we get
$$H^2-1+\frac{1}{a_0} \leq m_{C'} - 1 \leq H^2+1-\frac{1}{a_0}.$$
This implies $$H^2+\frac{1}{a_0} \leq m_{C'} \leq H^2 +2-\frac{1}{a_0}$$
 As $m_C$ is an integer, we therefore conclude
$$H^2< m_{C'}< H^2+2 \hspace{2mm} \text{and} \hspace{2mm} m_{C'} =  H^2+1 > 1.$$

Applying Riemann-Roch to $m_{C'}K_{X_n}+C'$ we get
\begin{align*}
& \chi(m_{C'}K_{X_n}+C') = h^0(m_{C'}K_{X_n}+C')-h^1(m_{C'}K_{X_n}+C')+h^2(m_{C'}K_{X_n}+C')\\
& = \frac{1}{2}\Bigl[\bigl(m_{C'}K_{X_n}+C'\bigr)\bigl(m_{C'}K_{X_n}+C'-K_{X_n}\bigr)\Bigr]+\chi(\mathcal{O}_X)\\
& = \frac{1}{2}\Bigl[\bigl(m_{C'}K_{X_n}+C'\bigr)\bigl((m_{C'}-1)K_{X_n}+C'\bigr)\Bigr] +\chi(\mathcal{O}_X)\\
& = \frac{1}{2}m_{C'}(m_{C'}-1)K_{X_n}^2+\frac{1}{2}m_{C'}K_{X_n}\cdot C'+\frac{1}{2}(m_{C'}-1)K_{X_n}\cdot C'+\frac{1}{2}C'^2+\chi(\mathcal{O}_X)\\
& = \frac{1}{2}m_{C'}(m_{C'}-1)K_{X_n}^2 +\frac{1}{2}\bigl(2m_{C'}-1\bigr)K_{X_n}\cdot \bigl(C+a_0H\bigr) + \frac{1}{2}\bigl(C+a_0H)^2+\chi(\mathcal{O}_X)\\
& =  \frac{1}{2}m_{C'}(m_{C'}-1)K_{X_n}^2 +\frac{1}{2}\bigl(2m_{C'}-1\bigr)K_{X_n}\cdot (a_0H) + \frac{1}{2}\bigl(2m_{C'}-1\bigr)K_{X_n}\cdot C\\
& \hspace{2cm} + \frac{1}{2}\bigl(C^2+2a_0C\cdot H+a_0^2H^2\bigr)+\chi(\mathcal{O}_X)\\
& = \frac{1}{2}H^2(H^2+1)K_{X_n}^2 - \frac{1}{2}a_0^2\bigl(2H^2+1\bigr)+\frac{1}{2}\bigl(2H^2+1\bigr)\bigl(2p_a(C)-2-C^2\bigr)\\
& \hspace{2cm} +\frac{1}{2}C^2  +a_0C\cdot H+\frac{a_0^2H^2}{2}+\chi(\mathcal{O}_X)\\
& = \frac{1}{2}H^2(H^2+1)K_{X_n}^2 - \frac{(2H^2+1)}{2}a_0^2 + \bigl(2H^2+1\bigr)\bigl(p_a(C)-1\bigr)  \\
& \hspace{2cm} - \frac{1}{2}C^2(2H^2+1) +\frac{1}{2}C^2+a_0C\cdot H+\frac{a_0^2H^2}{2}+\chi(\mathcal{O}_X)
\end{align*}
This implies
\begin{align*}
 \hspace{1cm} C^2 = & \frac{1}{2}(H^2+1)K_{X_n}^2 - \frac{(2H^2+1)}{2H^2}a_0^2 + \frac{\bigl(2H^2+1\bigr)}{H^2}\bigl(p_a(C)-1\bigr) + a_0\frac{C\cdot H}{H^2}+\frac{a_0^2}{2}\\
 &\hspace{2cm} +\frac{\chi(\mathcal{O}_X)}{H^2}
  - \frac{1}{H^2}\Bigl\{h^0(m_{C'}K_{X_n}+C')-h^1(m_{C'}K_{X_n}+C')+h^2(m_{C'}K_{X_n}+C')\Bigr\}.
\end{align*}
Note that by Lemma \ref{lem2.2} we have $h^0(m_{C'}K_{X_n}+C') = 0$ and $h^2(m_{C'}K_{X_n}+C') \leq 1.$
Therefore we have
\begin{align}\label{seq 7}
 C^2 & \geq \frac{1}{2}(H^2+1)K_{X_n}^2-a_0^2\frac{(1+H^2)}{2H^2}+\frac{(2H^2+1)}{H^2}\bigl(p_a(C)-1\bigr)+\frac{a_0C\cdot H +\chi(\mathcal{O}_X)-1}{H^2}\\
 &\geq  \frac{1}{2}(H^2+1)\bigl(K_{X}^2 -n\bigr) - a_0^2 -3 + \frac{a_0C\cdot H}{H^2}
\end{align}
\end{enumerate}
\vspace{2mm}

Combining Case 1 and Case 2 we have the result.
 \end{proof}

\end{thm}
\begin{thm}\label{thm4.5}
 Let $X$ be a smooth projective irreducible surface with $\chi(\mathcal{O}_X) < 1$  and  $-K_X$ is effective. Furthermore, let $L$ be a very ample divisor on $X$ with the intersection number $L\cdot K_X = - a_0 < 0$. Let $\pi_n : X_n = \Bl_n X \longrightarrow X$ be the blow-up map at $n$ points, and $H=\pi_n^*L$. Then for any integral curve $C$ on $X_n$, we have the following two cases :
 \begin{enumerate}
  \item \bf Case 1 : \it If $K_X^2\leq n$, then $$C^2 \geq \min\bigl\{\mathcal{L},\mathcal{M}\bigr\},$$
  \vspace{2mm}

  \item \bf Case 2 : \it If $K_X^2> n$, then
$$C^2 \geq \min\bigl\{\mathcal{M},\mathcal{N}\bigr\},$$
where
$$\mathcal{L}= \chi(\mathcal{O}_X)+\frac{(C\cdot H +a_0)}{2a_0}\Bigl(K_X^2-n\Bigr)-4, $$
$$\mathcal{N}=  \chi(\mathcal{O}_X)  + \frac{(C\cdot H+1)}{2a_0}\Bigl(K_X^2-n\Bigr) -4,$$
  and $$\mathcal{M}= \frac{1}{2}(H^2+1)\bigl(K_{X}^2 -n\bigr) - a_0^2 -3 + \frac{a_0C\cdot H+\chi(\mathcal{O}_X)-1}{H^2}.$$
  \vspace{1mm}

 \end{enumerate}
\begin{proof}
The proof is similar to the proof of Theorem \ref{thm3.2}. Let $C$ be an integral curve on $X_n$. We consider the following two cases :
\begin{enumerate}
 \item \bf Case 1 : \rm Let $m_C\neq 1$. Then procedding similarly as in Case 1 of Theorem \ref{thm3.2}  (see (\ref{seq1})) we have
 $$C^2 \geq \frac{\chi(\mathcal{O}_X)-1}{m_C-1} + \frac{1}{2}m_CK_{X_n}^2+2p_a(C)-3.$$
 This implies
 \begin{align}
  C^2\geq \chi(\mathcal{O}_X)+\frac{1}{2}m_CK_{X_n}^2-4.
 \end{align}
Now we consider the following two subcases :
\begin{enumerate}
 \item \bf Subcase 1 : \rm Suppose $K_X^2 \leq n$. Then, in this case, using Lemma \ref{lem2.2} we conclude
 \begin{align}
  C^2\geq \chi(\mathcal{O}_X)+\frac{(C\cdot H +a_0)}{2a_0}\Bigl(K_X^2-n\Bigr)-4
 \end{align}
\item \bf Subcase 2 : \rm Suppose $K_X^2 > n$. Then in this case, using Lemma \ref{lem2.2} we have
\begin{align}
 C^2 \geq \chi(\mathcal{O}_X)  + \frac{(C\cdot H+1)}{2a_0}\Bigl(K_X^2-n\Bigr)-4
\end{align}
\end{enumerate}
\item \bf Case 2 : \rm Let $m_C = 1$. Consider the  curve $C'=C+a_0H$ and proceeding similarly as in Case 2 of Theorem \ref{thm3.2} (see the inequality $(\ref{seq 7}))$, we get the following
\vspace{1mm}

$$C^2\geq \frac{1}{2}(H^2+1)\bigl(K_{X}^2 -n\bigr) - a_0^2 -3 + \frac{a_0C\cdot H+\chi(\mathcal{O}_X)-1}{H^2}.$$
\vspace{1mm}

We omit the details in Case 2. Now combining Case 1 and Case 2, we have the desired result.
\end{enumerate}
\end{proof}
\end{thm}
\begin{exm}\label{exm3.6}
Let $\pi_n : X_n = \Bl_n(\mathbb{P}^2) \longrightarrow \mathbb{P}^2$ be the blow-up of the projective plane $\mathbb{P}^2$ at $n$ points. We fix the notation $H$ for the pullback of a line $L$ in $\mathbb{P}^2$ under the map $\pi_n$. Then \begin{center}
$K_{X_n} = -3H+E_1+E_2+\cdots+E_n,$
and $a_0 = -K_{\mathbb{P}^2}\cdot L = 3.$
\end{center}
 Also, note that $-K_{\mathbb{P}^2} = \mathcal{O}_{\mathbb{P}^2}(3)$ is effective, and $\chi(\mathcal{O}_{\mathbb{P}^2}) =1$. Therefore, using Theorem \ref{thm3.2}, for any integral curve $C$ in $X_n$, we have the following inequalities:
\begin{itemize}
 \item \bf Case 1 : \it Suppose $K_{\mathbb{P}^2}^2 - n = 9 - n \leq 0.$
 Then $$C^2 \geq \min\Bigl\{ \Bigl(\frac{C\cdot H+3}{6}\Bigr)\bigl(9-n\bigr) - 3, \hspace{2mm} -n-3+3C\cdot H \Bigr\}.$$

 \item \bf Case 2 : \it Suppose $K_{\mathbb{P}^2}^2 - n = 9-n > 0$.
 Then $$C^2 \geq \min\Bigl\{ \Bigl(\frac{C\cdot H+1}{6}\Bigr)\bigl(9-n\bigr)-3,\hspace{2mm} -n-3+3C\cdot H \Bigr\}.$$
\end{itemize}
\end{exm}
\begin{exm}\label{exm3.7}
Let $\pi : \mathbb{F}_{e} = \mathbb{P}_{\mathbb{P}^1}\bigl(\mathcal{O}_{\mathbb{P}^1}\oplus \mathcal{O}_{\mathbb{P}^1}(-e)\bigr)\longrightarrow \mathbb{P}^1$ be a Hirzebruch surface for some non-negative integer $e$. Consider the blow-up map $\pi_n : \mathbb{F}_{e,n} := \Bl_n\mathbb{F}_e\longrightarrow \mathbb{F}_e$ at $n$ points. Let $C_0$ be the numerical class of the tautological line bundle  $\mathcal{O}_{\mathbb{F}_e}(1)$ and $f$ denotes the numerical class of a fiber of the map $f$. Then
$$K_{\mathbb{F}_e} \equiv -2C_0+(-2-e)f \hspace{2mm} \text{and} \hspace{2mm} K_{\mathbb{F}_e}^2 =8.$$ Let us fix a very ample line bundle $L \equiv C_0+(e+1)f$ on $\mathbb{F}_e$. Then $$a_0 = - K_{\mathbb{F}_e}\cdot L = 4+e \hspace{2mm} \text{and} \hspace{2mm} H^2 := \pi_nL^2 = e+1.$$  Also, note that $-K_{\mathbb{F}_e} \equiv 2C_0+(e+2)$ is effective, and   $\chi(\mathcal{O}_{\mathbb{F}_e}) = 1.$
Therefore, using Theorem \ref{thm3.2}, for any integral curve $C$ in $\mathbb{F}_{e,n}$, we have the following inequalities :
\begin{itemize}
 \item \bf Case 1 : \it Suppose $K_{\mathbb{F}_e}^2 - n = 8 - n \leq 0.$
 Then $$C^2 \geq \min\Bigl\{\frac{(C\cdot H +4+e)}{2(4+e)}(8-n)-3, \hspace{2mm} \frac{1}{2}(e+1)(8-n)-(4+e)^2-3+\frac{(4+e)(C\cdot H)}{(e+1)}\Bigr\}.$$

 \item \bf Case 2 : \it Suppose $K_{\mathbb{F}_e}^2 - n = 8-n > 0$.
 Then $$C^2 \geq \min\Bigl\{\frac{(C\cdot H +1)}{2(4+e)}(8-n)-3, \hspace{2mm}\frac{1}{2}(e+1)(8-n)-(4+e)^2-3+\frac{(4+e)(C\cdot H)}{(e+1)}\Bigr\}.$$
\end{itemize}
\end{exm}
\begin{exm}\label{exm3.8}
Let $B$ be a smooth projective curve over the complex numbers of genus $g\geq 1$ and $E = \mathcal{O}_B\oplus \mathcal{L'}$ be a rank 2 vector bundle on $B$ with $\deg(\mathcal{L}') < 3-3g$. Let $\pi : X = \mathbb{P}_B(E)\longrightarrow B$ be the associated ruled surface on $B$. We denote the class of tautological line bundle $\mathcal{O}_{\mathbb{P}_B(E)}(1)$ by $C_0$ and the class of a fiber of $\pi$ by $f$. By abuse of notation we use $C_0+Df$ to denote the divisor $\mathcal{O}_X(1)\otimes \pi^*D$ on $X$.

We choose a very ample divisor $D$ on $B$ having degree $\deg(D) = 2g+1-\deg(\mathcal{L}') > 2g+1.$

We also have $D+\det(E)$ is  very ample divisor on $B$ with $$\deg(D+\det(E)) = 2g+1-\deg(\mathcal{L}') + \deg(\mathcal{L}') \geq 2g+1.$$

Note that for any closed point $P\in B$, the divisors $D-P$ and $D+\det(E)-P$ are non-special divisors on the smooth curve $B$. This is indeed true as $$\deg(D-P) = \deg(D)-1 = 2g+1-\deg(\mathcal{L}')-1 = 2g-\deg(\mathcal{L}') > 2g-2,$$
and $$\deg(D+\det(E)-P) = \deg(D)+\deg(\mathcal{L}')-1  = 2g+1-\deg(\mathcal{L}') + \deg(\mathcal{L}') -1 > 2g-2.$$
Thus by \cite[Ex.2.11, Chapter 5, Section 2]{H77} we have $C_0+Df$ is very ample on $X=\mathbb{P}_B(E)$.

Then
\begin{align*}
a_0 & = -K_X\cdot (C_0+Df)\\
& = \bigl(2C_0-(2g-2+\deg(\mathcal{L}')\bigr)f\bigr)\cdot (C_0+Df)\\
&= 2\deg(\mathcal{L}')-2g+2-\deg(\mathcal{L}')+2\deg(D) \\
&= 2\deg(\mathcal{L}')-2g+2-\deg(\mathcal{L}')+ 4g+2-2\deg(\mathcal{L}')\\
& = 2g+4-\deg(\mathcal{L}') > 0.
\end{align*}
Also $$H^2 = \bigl(C_0+Df\bigr)^2 = \deg(\mathcal{L}')+2\bigl(2g+1-\deg(\mathcal{L}')\bigr) = 4g+2-\deg(\mathcal{L}').$$

Next we will show that the anti-canonical line bundle $-K_{\mathbb{P}_B(E)} \equiv 2C_0+\pi^*(K_B^{\vee}+\det(E)^{\vee})$ is an effective divisor on $X$.  It is enough to prove that $$H^0\bigl(B, \Sym^2(E)\otimes \det(E)^{\vee} \otimes K_B^{\vee}\bigr) \neq 0.$$
Note that $\mathcal{L}'^{\vee}\otimes K_B^{\vee}$ is a subbundle of $\Sym^2(E)\otimes \det(E)^{\vee} \otimes K_B^{\vee}$.\\
Now by Riemann-Roch theorem for line bundles we have $$h^0(\mathcal{L}'^{\vee}\otimes
K_B^{\vee}) \geq \deg(\mathcal{L}'^{\vee}\otimes K_B^{\vee})+
1 - g \geq 2-2g - \deg(\mathcal{L}') +1 -g = 3-3g
- \deg(\mathcal{L}') > 0. $$
Thus we conclude $$H^0\bigl(B, \Sym^2(E)\otimes \det(E)^{\vee} \otimes K_B^{\vee}\bigr) \neq 0.$$
This completes the proof of our claim, i.e. $-K_{\mathbb{P}_B(E)}$ is effective.

We also observe that $\chi(\mathcal{O}_X) = 1- g \leq 0$ and $K_X^2 = 8(1-g) \leq 0 \leq n$. Thus by Case 1 of the previous Theorem \ref{thm4.5} we have for any integral curve $C$ on the blow-up $\pi_n : X_n\longrightarrow X$
$$C^2\geq \min\{\mathcal{L},\mathcal{M}\}$$ where
\begin{align*}
\mathcal{L} & = \chi(\mathcal{O}_X)+\frac{(C\cdot H +2g+4-\deg(\mathcal{L}'))}{2\bigl(2g+4-\deg(\mathcal{L}')\bigr)}\Bigl(8(1-g)-n\Bigr)-4\\
&= \frac{(C\cdot H +2g+4-\deg(\mathcal{L}'))}{2\bigl(2g+4-\deg(\mathcal{L}')\bigr)}\Bigl(8(1-g)-n\Bigr) -g-3
\end{align*}
and
\begin{align*}
 \mathcal{M} = &\frac{1}{2}\bigl(4g+3-\deg(\mathcal{L}')\bigr)\Bigl(8(1-g)-n\Bigr)-(2g+4-\deg(\mathcal{L}'))^2-4\\
 & \hspace{3cm} +\frac{\bigl(2g+4-\deg(\mathcal{L}')\bigr)C\cdot H + \chi(\mathcal{O}_X)-1}{4g+2-\deg(\mathcal{L}')}.
\end{align*}
\end{exm}

\section{Weak bounded negativity for families of surfaces}
Let $f : Y=\mathbb{P}_X(E)\longrightarrow X$ be a projective bundle corresponding to a vector bundle $E$ of rank 3  on $X$. Then $Y$ is a family of surfaces isomorphic to $\mathbb{P}^2_{\mathbb{C}}$. Let $C$ be an integral vertical curve $C$ (i.e. $f(C)= x\in X$). Then $$C^2\geq 1.$$
  (Here the self-intersection of $C$ taken in the fiber $Y_x$ over the point $x\in X$).

  The following question is asked in \cite[Problem 3.3.6]{BBCRDHJKKRR12}.
\begin{prob}
 Let $f : Y \longrightarrow B$ be a morphism from a smooth projective threefold $Y$ to a smooth curve $B$ such that the general fibre is a smooth surface. Is there a
constant $b(Y,g)$ such that
$$C^2\geq-b(Y,g)$$
for all vertical curves $C \subset Y$ (i.e., $f(C)$ = a point) of geometric genus $g$? (Here the
self-intersection is computed in the fibre of $f$ containing $C$.)
\end{prob}

The above mentioned question motivates us to prove Theorem \ref{thm5.4}. We begin with the following important lemma.
\begin{lem}\label{lem1.1}
 Let $X$ be a smooth projective surface $X$ over $\mathbb{C}$ with $h^0(-K_X)>0$ and let $C$ be any integral curve of $X$. Then
 $$C^2\geq\textrm{min}\Bigl\{-2,\hspace{1mm} \chi(\mathcal{O}_X)+K_X^2-h^0(-K_X)-3\Bigr\}$$
 \begin{proof}
  Let $-K_X=P+N$ be the integral Zariski decomposition of the effective integral divisor $-K_X$. Let $C$ be an integral curve on $X$ which is not a component of the negative part of $-K_X$, then $-K_X\cdot C\geq 0$. By the adjunction formula, we have $C^2=2p_a(C)-2+(-K_X\cdot C)$. Hence $C^2\geq -2$.

  Now let $C$ be a component of $N$ which is the negative part of $-K_X$. Let $D=-K_X-C$ and $D=P+(N-C)$ is the integral Zariski decomposition of $D$. Then $$h^0(D)= h^0(-K_X-C) = h^0(P)=h^0(-K_X)>0.$$ Note that $h^0(K_X)=0$ and $h^0(K_X+C)=0$, hence $h^0(2K_X+C)=0$.
  \vspace{2mm}

  Now we apply Riemann-Roch for the divisor $-K_X-C$ on the smooth surface $X$, and we get
  \begin{align*}
  &h^0(-K_X-C)-h^1(-K_X-C)+h^2(-K_X-C)\\ & = h^0(-K_X-C)-h^1(-K_X-C)+h^0(2K_X+C)\\
  & =\chi(\mathcal{O}_X)+\frac{(-K_X-C)(-K_X-C-K_X)}{2}.
  \end{align*}
  This implies
  \begin{align*}
  h^0(-K_X-C)-h^1(-K_X-C)& =\chi(\mathcal{O}_X)+\frac{(K_X+C)(2K_X+C)}{2}\\
  &=\chi(\mathcal{O}_X)+K_X^2+3p_a(C)-3-C^2
  \end{align*}
  Therefore finally we have
  \begin{align*}
  C^2&\geq \chi(\mathcal{O}_X)+K_X^2-3-h^0(-K_X).
  \end{align*}
  Hence the result follows. Note that $h^0(-K_X) = h^0(P) = h^0(-K_X-C).$
 \end{proof}
 \end{lem}

\begin{prop}
 Let $f:Y\rightarrow B$ be a morphism from a smooth projective threefold $Y$ to a smooth curve $B$ such that generic  fibers are smooth surface i.e. $Y$ is a  fiber space. Let $Y_x=f^{-1}(x)$ be smooth fibers of the map $f$. Then the numbers $\chi(\mathcal{O}_{Y_x})$, $K_{Y_x}^2$ and $c_2(Y_x)$  are independent of $x\in B$.
 \begin{proof}
First note that two fibres are numerically equivalent. Hence $c_1(Y_x)$ is independent of $x$. Now consider the short exact sequence
$$0\rightarrow \mathcal{O}(-Y_x)\rightarrow \mathcal{O}_Y\rightarrow \mathcal{O}_{Y_x}\rightarrow 0,$$
such that $\chi(\mathcal{O}_{Y_x})=\chi(\mathcal{O}_Y)-\chi(\mathcal{O}(-Y_x))$. By Hirzebruch Riemann-Roch, $$\chi(\mathcal{O}(-Y_x))=\sum\limits_{j=0}^n\text{ch}_{n-j}(\mathcal{O}(-Y_x)).\text{td}_j(X),$$
where $\text{ch}_{n-j}(\mathcal{O}(-Y_x))$ is the Chern character of $\mathcal{O}(-Y_x)$ which depends only on $c_1(-Y_x)$. Therefore $\chi(\mathcal{O}_{Y_x})$ is independent of $x$.

Note that $c_2(Y_x)=c_2(\mathcal{T}_{Y_x})$, where $\mathcal{T}_{Y_x}$ is the tangent bundle of $Y_x$. We have the following exact sequence
$$0\rightarrow \mathcal{T}_{Y_x}\rightarrow \mathcal{T}_{Y}\otimes\mathcal{O}_{Y_x}\rightarrow \mathcal{N}_{Y_x/X}\rightarrow 0,$$
where $\mathcal{T}_{Y}$ is the tangent bundle of $Y$ and ${N}_{Y_x/X}\simeq\mathcal{O}(Y_x)\otimes\mathcal{O}_{Y_x}$ is the corresponding normal bundle. Now $c_2(\mathcal{T}_{Y_x})$ is a function of $c_1(\mathcal{N}_{Y_x/X})$, $c_1(\mathcal{T}_{Y}\otimes\mathcal{O}_{Y_x})$ and $c_2(\mathcal{T}_{Y}\otimes\mathcal{O}_{Y_x})$ which are independent of $x$.

Finally, the intersection product $K_{Y_x}^2=K_Y\cdot K_Y\cdot Y_x$ is also independent of $x$.
 \end{proof}
\end{prop}

 \begin{prop}\label{prop5}
  Now let $f:Y\rightarrow B$ be a morphism from a smooth projective threefold $Y$ to a smooth curve $B$ such that generic  fibers are smooth surface. Then $h^0(-K_{Y_x})$ is bounded above for all smooth fiber $Y_x$, $x\in B$.
  \begin{proof}
  Let $K_Y$ be the canonical divisor of $Y$ and $K_{Y_x}$ be the canonical divisor of smooth fiber $Y_x$.  By adjunction formula we have $K_{Y_x}=(K_{Y})\vert_{Y_x}$ as well as $-K_{Y_x}=(-K_{Y})\vert_{Y_x}$.  Then by using semicontinuity  theorem (\cite{H77}, III, Theorem 12.8), we conclude that $h^0(Y_x,-K_{Y_x})$ is bounded above for all smooth fiber $Y_x$ and $x\in U$, where $U$ is some open subset of $B$. But $B\setminus U$ is a finite set. So eventually we conclude that $h^0(Y_x,-K_{Y_x})$ is bounded above for all smooth fiber $Y_x$.
\end{proof}
\end{prop}

 \begin{thm}\label{thm5.4}
  Let $f:Y\rightarrow B$ be a morphism from a smooth projective threefold $Y$ to a smooth curve $B$ such that generic  fibers are smooth surface. Let $Y_x=f^{-1}(x)$ be a smooth fiber for some $x\in B$, and $h^0(-K_{Y_x})\leq l$ for all $x\in B$ and for some positive integer $l$. Let $C\subset Y$ be any integral vertical curve (i.e. $f(C)=$ point) of geometric genus $p_g(C)$ such that $C$ is contained in a smooth fiber $Y_x$. Then $$C^2\geq b(Y,g)$$
  where $b(Y,g)= \textrm{min}\Bigl\{-2, \chi(\mathcal{O}_{Y_x})+K_{Y_x}^2-3-l,  K_{Y_x}^2+\chi(\mathcal{O}_{Y_x})-3,K_{Y_x}^2-3c_2(Y_x)+2-2p_g(C)\Bigl\}$.\\
  $($Here the self-intersection is computed in the fiber of $f$ containing $C.)$
  \begin{proof}
   Let   $K_{Y_x}$ be the canonical divisor of  a smooth fiber $Y_x$ containing the integral curve $C$.

   Consider the following two cases :
   \begin{itemize}
    \item \bf Case 1 : \rm If $h^0(Y_x,-K_{Y_x})> 0$, then applying Proposition \ref{prop5} and using Lemma \ref{lem1.1} we have $$C^2\geq  \textrm{min}\Bigl\{-2,\hspace{1mm} \chi(\mathcal{O}_{Y_x})+K_{Y_x}^2-3-l)\Bigr\}.$$
    \item \bf Case 2 : \rm If $h^0(Y_x,-K_{Y_x})=0$, then using Theorem \ref{thm4} we get $$C^2\geq \min \Bigl\{K_{Y_x}^2+\chi(\mathcal{O}_{Y_x})-3, \hspace{1mm} K_{Y_x}^2-3c_2(Y_x)+2-2p_g(C)\Bigr\}.$$
   \end{itemize}
   Combining both the cases, we have the desired inequality.
\end{proof}
 \end{thm}

 \begin{corl}\label{corl4.6}
Let $f:Y\rightarrow B$ be a morphism from a smooth projective threefold $Y$ to a smooth curve $B$ such that generic  fibers are smooth surfaces. Let $C=\sum\limits_{i=1}^m C_i\subset Y$ be any reduced vertical curve such that $C$ is contained in a smooth fiber $Y_x$. Let $g=\max\bigl\{p_g(C_i) \mid 1\leq i \leq m\bigr\}$. Then there exists a constant $B(Y,g)$ only depending on $Y$ and
$g$, such that
$$C^2 \geq  B(X,g)$$.
\begin{proof}
The proof is identical to the proof of \cite[Theorem 2.1]{H19}. We omit the details.
\end{proof}
 \end{corl}

\end{document}